\documentclass[fleqn,preprint,12pt,3p]{elsarticle}
\usepackage{amsmath,amssymb,amsthm}
\usepackage{lineno,hyperref}
\usepackage{graphicx}
\usepackage[normalem]{ulem}
\usepackage{amssymb,amsmath,array,amsthm}
\usepackage{subfig}
\newtheorem{theorem}{Theorem}[section]

\usepackage{color}
\usepackage{ulem}

\newtheorem{remark}[theorem]{Remark}

\usepackage{soul,xcolor}
\usepackage[]{algorithm2e}
\RestyleAlgo{boxruled}

\usepackage{color}
\begin{document}
\begin{frontmatter}
\title{Local recovery bounds for prior support constrained Compressed Sensing}
\author{K. Z. Najiya, Munnu Sonkar and C. S. Sastry \\
Department of Mathematics \\
Indian Institute of Technology, Hyderabad, 502285, India. \\ Email:\{ma17resch11004,  ma17resch01001 and csastry\}@iith.ac.in}
\date{}
\begin{abstract}
Prior support constrained compressed sensing has of late become popular due to its potential for applications. The existing results on recovery guarantees provide global recovery bounds in the sense that they deal with full support. However, in some applications, one might be interested in the recovery guarantees limited to the given prior support, such bounds may be termed as local recovery bounds. The present work proposes the local recovery guarantees and analyzes the conditions on associated parameters that make recovery error small.
\end{abstract} 
\end{frontmatter}
\section{Introduction}
In Compressed Sensing (CS), a sparse signal $x \in \mathbb{R}^{n}$ can be recovered from a small set of linear measurements $y$ satisfying $y = Ax$ with $k \leq n$ where $k$ is the number of nonzero elements in $x$. The results guaranteeing the recovery performance  depend on the coherence, Restricted Isometry Property (RIP) and Restricted Orthogonality Constant (ROC) of the sensing matrix $A$ \cite{candes2005decoding}\cite{Kashin}\cite{EulerSq}. In many applications, one obtains some a priori information about the partial support of the sparse solution to be recovered. For instance, in interior reconstruction in Computed Tomography, one has before hand some apriori information corresponding to the support of the interior region \cite{Klann2015WaveletMF},\cite{localCBP}.
There are other applications, like recovering time-correlated signals \cite{vaswani2010modified}, wherein prior-support constrained sparse recovery attains importance. Of late, the support constrained CS  
has caught the attention of several researchers \cite{friedlander2011recovering}\cite{Klann2015WaveletMF}\cite{vaswani2010modified} to name a few. In \cite{vaswani2010modified}, the authors have modified the 1-norm function by taking zero weights on the known partial support, minimizing thereby the terms in the complement of prior support set. While in \cite{friedlander2011recovering}, by considering general values for weights, the authors have provided stability and robustness of weighted $l_1$-norm minimization problem in terms of Restricted Isometry Constant (RIC). The authors of \cite{liu2014compressed} have studied similar performance guarantees in terms of mutual coherence of associated sensing matrix.
The work embedded in \cite{Chen2016RecoveryOS} has provided a less restrictive sufficient condition and a tighter error bounds under some conditions for the weighted-$l_1$ norm problem in terms of RIC and ROC. More recently, the authors of \cite{Ge2018RecoveryOS} have provided the stability and robustness of the weighted-$l_1$ norm in terms of block RIP conditions on the underlying measurement matrix. 
\subsection{Motivation for our work}
\par In applications like interior tomography, however, one is interested in the recovery guarantees limited to the interior portion that is accounted for by partial support. This is because the recovery on interior portion only attains importance in such an application and the reconstruction outside interior portion  is, in general, bad \cite{Farrokhi}\cite{Klann2015WaveletMF}.  Motivated by this, 
%the need for error bound limited to the prior support set and uniqueness of solution in a general setting, 
the present work proposes a new local recovery bound, in the sense that it pertains only to the prior support. It is to be emphasized here that by a prior support set $T$, we mean an arbitrary subset of ``full support" $\{1,2,\ldots,n\}$, which, in general, does not have to  be fully contained in the true support of $x$. This is in contrast to the existing results that provide global bounds (that is, on entire solution support). Further, both analytically and empirically we demonstrate the conditions on associated parameters that reduce the reconstruction error. 
\par The paper is organized as 5 sections. In section 2, we provide basic introduction to Compressed Sensing, existing recovery bounds and a summary of our contribution. Section 3 presents the local recovery bound followed by an analysis and comparison in section 4. The paper ends with the concluding remarks in section 5.
\section{Compressed sensing}
Compressed sensing(CS) is a technique that reconstructs  a signal, which is compressible or sparse in some domain, from a small set of linear measurements.  Let $\sum_{k}^{n} := \{x \in \mathbf{R}^{n}:\|x\|_0 \leq k \}$ be the set of all k-sparse signals in $\mathbf{R}^{n}$. Here $\|x\|_0$ stands for the number of nonzero components in $x$. The best $k$-term approximation $x_k\in \mathbb{R}^n$ of $x$ retains at most $k$ largest magnitude coordinates of $x$, the rest of the coordinates are set to zero. For simplicity,  we denote $\{1,2, \dots, n\}$ by $[n]$.  For $A  \in \mathbf{R}^{m \times n}$ ( $m << n$) and  an error vector $\xi\in \mathbf{R}^m$,  suppose  $y=Ax+\xi$ such that $\|\xi\|_2 \leq \epsilon$. One may recover the sparsest solution of the noisy matrix system from the following minimization problem \cite{cai2010stable}:
\begin{equation}\label{l0}
   (P_0) \hspace{.4cm}  min \|x\|_0 \hspace{0.3cm} subject\hspace{0.3cm} to \hspace{0.3cm} (y-Ax) \in \mathcal{B},
\end{equation}
  where $\mathcal{B}$ is a bounded subset of $\mathbf{R}^m.$ For noiseless case,  $\mathcal{B} =\{0\}$ and for noisy case ${\mathcal{B}=\{\xi\in \mathbf{R}^m: \|\xi\|_2\leq \epsilon\}}$. Since $l_0$ minimization problem becomes  NP-hard as the dimension increases, the convex relaxation of $l_0$ problem (\ref{l0}) has been proposed as
\begin{equation}\label{eq:P1}
   (P_1) \hspace{.4cm}  min \|x\|_1 \hspace{0.3cm} subject\hspace{0.3cm} to \hspace{0.3cm} (y-Ax) \in \mathcal{B}.
\end{equation}
 The coherence $\mu(A)$ of a matrix $A$ is the largest absolute normalized inner product between different columns of it, that is, 
\begin{equation*}
    \mu(A)=  \max_{1\leq i,j \leq M, \; i \neq j}\frac{  |a_i^T a_j|}{\|a_i\|_2 \|a_j\|_2}, 
\end{equation*}
where $a_i$ denotes the $i^{th}$ column in $A$. 
For a $k$ sparse vector $x$, it is known \cite{elad2010sparse} that the following inequality holds:
\begin{equation} \label{MIP}
(1-(k-1)\mu)\|x\|_{{2}}^{2}\leq \|Ax\|_{{2}}^{2}\leq (1+(k-1)\mu )\|x\|_{{2}}^{2}.
\end{equation}
The $k$-th Restricted Isometry Constant ($k$-RIC)  $\delta_k$ of a matrix $A$ is the smallest number $\delta_k \in(0,1) $ such that 
\begin{equation*}
(1-\delta _{k})\|x\|_{{2}}^{2}\leq \|Ax\|_{{2}}^{2}\leq (1+\delta _{k})\|x\|_{{2}}^{2}, 
\end{equation*}
for all $k$-sparse vectors $x$. The Restricted Orthogonality Constant(ROC) $\theta_{s,\Tilde{s}}$ of a matrix $A$ is the smallest real number such that
$$
    |\eta'A_{T}'A_{\tilde T} \tilde \eta| \leq \theta_{s,\tilde s} \|\eta\|_2 \| \tilde \eta \|_2,  
$$
for all disjoint sets $T$ and $\tilde T$ with $|T| \leq s$ and $| \tilde T |\leq \tilde s $ such that $s+ \tilde{s}\leq n$ and for all vectors $\eta \in \mathbb{R}^{|T|}$ and $\tilde \eta \in \mathbb{R}^{|\tilde T|}$. Here, $A_T$ denotes the 
sub matrix of columns of $A$ restricted to the indices in $T \subseteq [n]$.
 D. Donoho and X. Huo \cite{donoho2001uncertainty} have  shown exact recovery condition in noiseless case for $P_1$ in terms of mutual coherence. 
If $x$ is $k$ sparse vector and matrix $A$ is $k$-RIP compliant, $ k<  \frac{1}{2}\big( 1+\frac{1}{\mu}\big)$ is an exact recovery condition for $P_1$ problem. Then T. Cai et. al. \cite{cai2010stable} have extended this result to noisy case.
\begin{theorem}(T. Cai et. al. \cite{cai2010stable}):
\label{CaiBound}
Consider the model $y=Ax+\xi$ with $\|\xi\|_2\leq \epsilon.$ Suppose $x$ is in $R^n$, $A\in R^{m \times n}$ and $x_k$ represents its best $k-term$ approximation with
\begin{equation} \label{boundK_classical}
    k<  \frac{1}{2} \biggl( 1+\frac{1}{\mu}\biggl).
\end{equation}
Let $x^*$ be the minimizer of $P_1$. Then $x^*$ obeys 
    \begin{equation*}
        \|x^* - x\|_2\leq C_0^{(0)} \epsilon + C_1^{(0)} \|x-x_k\|_1, 
    \end{equation*}
\text{where} $C_0^{(0)}=\frac{ 2\big(1+\mu - 2 \mu k+2\sqrt{\mu k (1+( k-1)\mu)}\big)}{(1+\mu - 2 \mu k)\sqrt{1+\mu}}$ and $
C_1^{(0)}=\frac{2(1+\mu) \sqrt{\mu}}{(1+\mu - 2 \mu k)\sqrt{1+\mu}}.$
\end{theorem}

\subsection{Compressed sensing with prior support constraint}
It may be noted that the reconstruction method given by $P_1$ in (\ref{eq:P1}) is non adaptive as no  information about $x$ is used in $P_1$. It can, however, be made partially adaptive by imposing constraints on the support of the solution to be obtained.  In                     \cite{friedlander2011recovering} \cite{liu2014compressed} \cite{vaswani2010modified}  (and the references therein) the authors have modified the cost function of $P_1$ problem by incorporating the prior support information into the reconstruction process as detailed in the following subsection.
\subsection{Previous Work}
\label{prev_work}
\par Consider that $T$ is the known partial support information of signal $x$, which is expected in recovered solution. Suppose $T_0$ stands for the support of the best $k$-term approximation $x_k$ of $x$, where $x$ is the actual solution of $y=Ax$.
 In \cite{vaswani2010modified}, the authors have modified the $P_1$ problem by considering zero weights in $T$ and have posed it as follows: 
\begin{equation}  \label{Nama}
min \|x_{T^c}\|_1 \; \text{ subject to} \; y=Ax.
\end{equation}
In the above problem, the weights are set to 1 on $T^c$ and to 0 on $T$. In \cite{friedlander2011recovering}, nevertheless, the authors have posed this problem for  a general weight vector and an arbitrary subset $T$ of $[n]$ the following way:
\begin{equation} \label{weighted}
    (P_{1,w})\hspace{.3cm} min \|x\|_{1,w} \hspace{0.3cm} subject\hspace{0.3cm} to \hspace{0.3cm} \|y-Ax\|\leq \epsilon,
\end{equation}
where  $\|x\|_{1,w}:=\sum_{i}w_i|x_i|$ with $w_i=w$ for $i \in T$ and $w=1$ for $i \notin T$,  $T$  can be drawn  from the estimate of the support of the signal $x$ or from its largest coefficients. The stability result proposed in \cite{friedlander2011recovering} is as follows: 
\begin{theorem}(M. Friedlander et. al.  \cite{friedlander2011recovering}): \label{Fried}
Let $x\in \mathbf{R}^n$ and let $x_k$ be its best $k-$term approximation, supported on $T_0$. Let $T \subset [n]$ be an arbitrary set and define $\rho, \alpha$ such that $|T|=\rho k$ and $|T_0 \cap T|=\alpha\rho k$. Suppose that there exists an $a\in \frac{1}{k} \mathbb{Z}$, with $a \geq (1-\alpha)\rho$, $a>1$, and the measurement matrix $A$ has RIP with 
\begin{equation}
    \delta_{ak}+\frac{a}{\beta^2}\delta_{(a+1)k}< \frac{a}{\beta^2}-1, 
\end{equation}
where $\beta=w+(1-w)\sqrt{1+\rho-2\alpha \rho}$ for some given $0\leq w\leq 1$. Then the solution $x^*$ to the (\ref{weighted}) obeys
\begin{equation} \label{FriedBound}
    \|x^* - x\|_2\leq C_0^{(1)} \epsilon + C_1^{(1)}  \big(w\|x-x_k\|_1 + (1-w)\|x_{T^c\cap T_0^{c}}\|_1\big),
\end{equation}

\noindent where
 \begin{equation}\label{coeff_frd}
\begin{split}
 C_0^{(1)} &=\frac{2\big(1+\frac{w+(1-w\sqrt{1+\rho-2 \alpha \rho})}{\sqrt{a}}\big)}{\sqrt{1-\delta_{(a+1)k}}-\frac{w+(1-w)\sqrt{1+\rho-2 \alpha \rho}}{\sqrt{a}}\sqrt{1+\delta_{ak}}}   \text{ and } \\
C_1^{(1)} & =\frac{2(ak)^{-1/2}\big(\sqrt{1-\delta_{(a+1)k}}+ \sqrt{1+\delta_{ak}}\;\big)}{\sqrt{1-\delta_{(a+1)k}}-\frac{w+(1-w)\sqrt{1+\rho-2 \alpha \rho}}{\sqrt{a}}\sqrt{1+\delta_{ak}}}.
\end{split}
\end{equation}
\end{theorem}
\qed \\
It has been shown in \cite{friedlander2011recovering} that a signal $x$ can be stably and robustly recovered from $P_{1,w}$ problem if  at least $50\%$ of the partial support information is accurate. 
It is worth a mention here that the above stability result has been proposed in terms of RIC-$\delta_k$. In \cite{liu2014compressed}, however, the authors have proposed a similar stability result, albeit in terms of coherence parameter, which is summarized as follows:
\begin{theorem}(Haixiao et. al.  \cite{liu2014compressed}): \label{IETbound}
Let $x \in \mathbf{R}^n$, and let $x_k$ be its best $k$-term approximation, supported on $T_0$. Let $T \subset [n]$ be an arbitrary set and define $\rho, \alpha$ such that $|T|=\rho k$ and $|T_0 \cap T|= \alpha\rho k$. Suppose that 
\begin{equation} \label{boundK_iet}
     k<  \frac{L}{2} \big( 1+\frac{1}{\mu} \big), 
\end{equation}
where $0\leq w\leq 1$, $L=\frac{Q+2-\sqrt{Q(Q+4)}}{1+w}$ and $Q=\frac{(1-w)^2(1+\rho-2\alpha\rho)}{1+w}$. Then the solution $x^*$ to (\ref{weighted}) obeys 
\begin{equation} \label{IETBound}
    \|x^* - x\|_2\leq C_0^{(2)} \epsilon + C_1^{(2)} \big(w\|x-x_k\|_1 + (1-w)\|x_{T^c\cap T_0^{c}}\|_1\big), 
\end{equation}
\end{theorem}
\noindent where 
\begin{equation}\label{coeff_iet}
 \begin{split}
 C_0^{(2)} &= \frac{2 \big(1-( k-1)\mu-\mu k w+(1+w)\sqrt{\mu k (1+(k-1)\mu)} \big )}{(1-(k-1)\mu -\mu k w)\sqrt{1+\mu} - \sqrt{\mu}(1+\mu)(1-w)\sqrt{\mu k(1+\rho-2\alpha \rho)}} \text{ and } \\
C_1^{(2)} & =\frac{2(1+\mu)\sqrt{\mu}}{(1-(k-1)\mu -\mu k w)\sqrt{1+\mu} - \sqrt{\mu}(1+\mu)(1-w)\sqrt{\mu k(1+\rho-2\alpha \rho)}}.
\end{split}
\end{equation}
\qed
\par The authors of \cite{Chen2016RecoveryOS} have proposed less restrictive sufficient conditions and a tighter bound with respect to the standard 1-norm problem under some conditions for the weighted 1-norm problem in terms of RIC and ROC. The stability result is as follows:
\begin{theorem}(Chen et. al.  \cite{Chen2016RecoveryOS})\label{ric_roc}
Let $x\in\mathbb{R}^n$ be an arbitrary signal and its best $k$-term approximation support on $T_0\subseteq [n]$ with $|T_0| \leq k$. Let $T\subseteq [n]$ be an arbitrary set and denote $\rho \geq 0$ and $0 \leq \alpha \leq 1$ such that $|T| =\rho k$ and $|T\cap T_0| =\alpha \rho k$. Let $y=Ax+z$ with $\|z\|_2\leq \epsilon$ and $x^*$ is the minimizer of (\ref{weighted}). If $\delta_a+C_{a,b,k}^{\alpha,w} \theta_{a,b} <1$ for some $a,b\in \mathbb{N}$ with $1\leq a \leq k$, where $C_{a,b,k}^{\alpha,w}=max\bigg\{ \frac{s}{\sqrt{ab}},\sqrt{\frac{s}{a}}\bigg\}$ with $s=k-a+wk+(1-w)\sqrt{(1+\rho-2\alpha\rho)k}.max\{\sqrt{(1+\rho-2\alpha\rho)k},\sqrt{a}\}$. Let $d=k$, for $w=1$ and $d$=max$\{k,(1+\rho-2\alpha\rho)k\}$, for $0\leq w <1$. Then
\begin{equation} \label{ric_roc_Bound}
    \|x^* - x\|_2\leq C_0^{(3)} \epsilon + C_1^{(3)}  \big(w\|x-x_k\|_1 + (1-w)\|x_{T^c\cap T_0^{c}}\|_1\big), 
\end{equation}

\noindent where
\begin{equation}\label{coeff_ric_roc}
   C_0^{(3)}=\frac{2\sqrt{2(1+\delta_a)d/a}}{1-\delta_a-C_{a,b,k}^{\alpha,w}\theta_{a,b}} ;\; C_1^{(3)}=\frac{2\sqrt{2d}C_{a,b,k}^{\alpha,w}\theta_{a,b}}{(1-\delta_a-C_{a,b,k}^{\alpha,w}\theta_{a,b})s}+\frac{2}{\sqrt{d}}. 
\end{equation}
\end{theorem}\qed

\noindent A vector $x=(\underbrace{x_1, \ldots,x_{d_1}}_\text{x[i]}, \ldots,\underbrace{x_{n-{d_M+1}}, \ldots,x_n}_\text{x[M]}\}\in\mathbb{R}^n$, where $x[i]$ is the $i^{th}$ block of $x$ of size $d_i$ w.r.t $\mathbb{T}=\{d_1,d_2, \ldots,d_M\}$ with $\sum_{i=1}^M d_i=n$, is said to be block $k$- sparse over $\mathbb{T}$ if the number of non-zero blocks in $x$ is at most $k$. Recently the authors of \cite{Ge2018RecoveryOS} have  introduced the following weighted block $l_{2}/l_{1}$-norm problem for a given $L$ disjoint prior block support estimates $T_j\subseteq [M]$ for $j=1,2,\ldots,L$ with $\cup_{j=1}^L T_i=T$ as 
\begin{equation}{\label{block_P1w}}
     \min_{x\in\mathbb{R}^n} \;\|x_{\bf{w}}\|_{2,1}=\sum_{i=1}^M  w_i\|x[i]\|_2 \hspace{0.3cm} subject\hspace{0.3cm} to \hspace{0.3cm} \|y-Ax\| \leq \epsilon,
\end{equation}
where ${\bf{w}}\in [0,1]^M$ is defined by ${\bf{w}}_i=1,\; i \in T^c$ and ${\bf{w}}_i=w_j; \; i \in T_j$ for $i\in [M]$.
 Note that when $d_i=1$ for all $i=1,\ldots,M$, and block sparsity reduces to the standard sparsity and if the number of support estimates $L=1$, then the weighted block $l_{2}/l_{1}$-norm problem (\ref{block_P1w}) reduces to the weighted 1-norm problem in (\ref{weighted}). In this particular case the stable recovery result of (\ref{block_P1w}) in \cite{Ge2018RecoveryOS} deduces to the following result:

\begin{theorem}(Ge et.al. \cite{Ge2018RecoveryOS})\label{block_thm}
For an arbitrary signal $x\in \mathbb{R}^n$, which satisfies $y=Ax+z$ with $\|z\|_2\leq \epsilon$, let $x_k$ be its best $k$-term approximation and $T_0=supp(x_k)$. Suppose that $x^*$ is the minimizer of (\ref{weighted}) and $T$ is the prior block support of $x$  satisfying $|T|= \rho k$,\;$|T \cap T_0| = \alpha\rho k$.
 If A satisfies the RIP with $\delta_{tk} < \sqrt{(t-d)/(t-d+\Upsilon_L^2)}$  for $t>d$, where,
$\Upsilon_L=w+(1-w)\sqrt{1+ \rho-2 \alpha\rho}$ and $d=1$ for $w=1$. For $0\leq w<1$, d=1 if $\alpha\geq1/2$ and $d=1+ \rho-2 \alpha\rho$ if $\alpha<1/2$. Then
\begin{equation} \label{block_Bound}
    \|x^* - x\|_2\leq C_0^{(4)} \epsilon + C_1^{(4)} \big(w\|x-x_k\|_1 + (1-w)\|x_{T^c\cap T_0^{c}}\|_1\big), 
\end{equation}
where 
\begin{equation}\label{coeff_block}
\begin{split}
    C_0^{(4)} &=\frac{2\sqrt{2(t-d)(t-d+\Upsilon_L^2)(1+\delta_{tk})}}{(t-d+\Upsilon_L^2)\bigg(\sqrt{\frac{t-d}{(t-d+\Upsilon_L^2)}}-\delta_{tk}\bigg)}  
    \\
    C_1^{(4)} &=\frac{2}{\sqrt{k}} \Bigg(
    \frac{\sqrt{2}\delta_{tk}\Upsilon_L+\sqrt{(t-d+\Upsilon_L^2)\bigg(\frac{t-d}{(t-d+\Upsilon_L^2)}-\delta_{tk}\bigg)\delta_{tk}}}{(t-d+\Upsilon_L^2)\bigg(\sqrt{(t-d+\Upsilon_L^2)\frac{t-d}{(t-d+\Upsilon_L^2)}}-\delta_{tk}\bigg)}+\frac{1}{\sqrt{d}}\Bigg).
    %X= \sqrt{2}\delta_{tk}\Upsilon_L+\sqrt{(t-d+\Upsilon_L^2)\bigg(\frac{t-d}{(t-d+\Upsilon_L^2)}-\delta_{tk}\bigg)\delta_{tk}}
\end{split}
\end{equation}
\end{theorem}

\section{Local recovery bounds} \label{Present-Work}
As discussed already, the present work deals with obtaining a recovery bound on $\|x^*_T - x_{T}\|$ for $T\subset[n]$. It is shown in the later part that when $T$ contains the indices corresponding to the largest magnitude entries in $x_k$, our error bound is much smaller than the global error bounds (\ref{FriedBound}),(\ref{IETBound}),(\ref{ric_roc_Bound}) and (\ref{block_Bound}) for weighted-1-norm problem.
Further, in most of the cases the sufficient condition on $k$ in this bound can be shown to be less pessimistic than the corresponding ones in the standard 1-norm (\ref{boundK_classical}) and weighted 1-norm cases (\ref{boundK_iet}).  Our contribution may be summarized as the following result: 

\begin{theorem} \label{theorem_hT}
Let $x$ be in $\mathbb{R}^n$ satisfying $\|y-Ax\|_2\leq \epsilon$ where $A\in \mathbb{R}^{m\times n}$ and $y\in \mathbb{R}^{m}$ with $m<n$. Let $x_k$ be its $k$-term approximation  supported on $T_0$. Let $T \subset [n]$ be an arbitrary set. Define $\rho$ and $\alpha$ such that $|T|=\rho k$ and $|T \cap T_0|=\alpha \rho k$. Suppose that 
\begin{equation}\label{boundK_new}
k< 
\begin{cases}
\Big(\frac{1}{2 \sqrt{\rho} (2 w \sqrt{\alpha}+1)}( w+\sqrt{w^2+4(2w \sqrt{\alpha}+1)(1+\frac{1}{\mu}) })\Big)^2,
\text{~if~} w \in (0, 1] \\
\frac{1}{\rho}\Big(1+\frac{1}{\mu}\Big), 
\text{~if~} w=0,
\end{cases}
\end{equation}
\noindent then the solution $x^*$  on $T$ to 
$(\ref{weighted})$ obeys 
\begin{equation}\label{error bound T}
    \|x^*_{T}-x_{T}\|_2 \leq C_0 \epsilon+ C_1\bigl( w\|x-x_k\|_1 +(1-w)\|x_{T^c \cap T_0^c}\|_1+\|x_{T^c \cap T_0}\|_1 \bigl),
\end{equation}
\noindent where
\begin{equation}\label{coeff}
C_0=\frac{2 \sqrt{1+(\rho k-1)\mu}}{(1+\mu+w\mu\sqrt{\rho k}-\mu \rho k(2w\sqrt{\alpha}+1))}, 
C_1=\frac{2\mu\sqrt{\rho k}}{(1+\mu+w \mu\sqrt{\rho k}-\mu \rho k(2w\sqrt{\alpha}+1))}.
\end{equation}
\end{theorem}

\noindent \textbf{Proof:}
Suppose $h=x^*-x$. From the definition of $x^*$, we have
\begin{equation*}
   \begin{split}
       \|{ x+h}\|_{1,w} & \leq \|x\|_{1,w}\\
       w\|x_{T\cap T_0}+h_{T\cap T_0}\|_1+w\|x_{T\cap T_0^c} +h_{T\cap T_0^c}\|_1 &\\ +\|x_{T^c}+h_{T^c}\|_1 & \leq 
w\|x_{T\cap T_0}\|_1+w\|x_{T\cap T_0^c}\|_1+\|x_{T^c}\|_1 \\
w\|x_{T\cap T_0}\|_1-w \|h_{T\cap T_0}\|_1+w \|h_{T\cap T_0^c}\|_1 &\\ -w \|x_{T\cap T_0^c}\|_1  +\|h_{T^c}\|_1-\|x_{T^c}\|_1 & \leq 
w\|x_{T\cap T_0}\|_1+\|x_{T\cap T_0^c}\|_1  +\|x_{T^c}\|_1\\
\|h_{T^c}\|_1 & \leq w \|h_{T\cap T_0}\|_1- w ||h_{T\cap T_0^c}||_1 \\&  +2(w \|x_{T\cap T_0^c}\|_1 +\|x_{T^c}\|_1)\|_1
   \end{split} 
\end{equation*}
Consider that  $e= w \|x_{T\cap T_0^c}\|_1 +\|x_{T^c}\|_1$. Then, we have 
\begin{equation} \label{h_tc}
    \|h_{T^c}\|_1 \leq w \|h_{T\cap T_0}\|_1- w ||h_{T\cap T_0^c}||_1+2e.
\end{equation}
Now, since $\langle Ah_{T^c},Ah_{T}\rangle  =\langle Ah,Ah_{T}\rangle - \langle Ah_{T},Ah_{T}\rangle$ and by the inequality (\ref{MIP}), it follows that
\begin{equation}\label{eqqq2}   \begin{split} 
     (1-(\rho k-1)\mu) \|h_{T}\|_2^2 \leq  \mid \langle Ah_{T},Ah_{T}\rangle\mid
     & \leq \mid\langle Ah,Ah_{T}\rangle\mid+\mid \langle Ah_{T^c},Ah_{T}\rangle\mid\\
     & \leq \|Ah\|_2\|Ah_{T}\|_2+\sum_{i\in T^c,j\in T} \mid\langle A_i,A_j\rangle\mid \; \mid h_i h_j \mid\\
     & \leq 2 \epsilon \sqrt{1+(\rho k-1)\mu}\|h_{T}\|_2+\mu\|h_{T^c}\|_1
     \|h_{T}\|_1.
    \end{split}
\end{equation}
In view of the inequalities $||h_{T\cap T_0}||_1\leq \sqrt{|T\cap T_0|}  ||h_{T\cap T_0}||_2 = \sqrt{\alpha \rho k } ||h_{T}||_2$, $||h_{T}||_2\leq ||h_{T}||_1$ and the ones in  (\ref{h_tc}) and (\ref{eqqq2}),  
we have the following:
\begin{equation*}
\begin{split}
    (1-(\rho k-1)\mu) \|h_{T}\|_2
& \leq 2 \epsilon \sqrt{1+(\rho k-1)\mu}+\mu\sqrt{\rho k} \|h_{T^c}\|_1 \\ 
& \leq 2 \epsilon \sqrt{1+(\rho k-1)\mu} +\mu\sqrt{\rho k} (w\|h_{T\cap T_0}\|_1-w\|h_{T\cap T_0^c}\|_1+2e)\\
 & \leq 2 \epsilon \sqrt{1+(\rho k-1)\mu} +\mu\sqrt{\rho k} (2w \|h_{T\cap T_0}\|_1-w\|h_{T}\|_1+2e)\\
 %& \leq 2 \epsilon \sqrt{1+(\rho k-1)\mu}+\mu\sqrt{\rho k} (2w \sqrt{\alpha \rho k} \|h_{T\cap T_0}\|_2-w \|h_{T}\|_2+2e)\\
 &\leq 2 \epsilon \sqrt{1+(\rho k-1)\mu}+2\mu\sqrt{\rho k}e+\mu\sqrt{\rho k} (2 w \sqrt{\alpha \rho k} -w )\|h_{T}\|_2, %\\
\end{split}
\end{equation*}
which results in 
\begin{equation*}
    (1+\mu+w \mu\sqrt{\rho k}-\mu \rho k(2w \sqrt{\alpha}+1)) \|h_{T}\|_2 \leq  2 \epsilon \sqrt{1+(\rho k-1)\mu} +2\mu\sqrt{\rho k}e.
\end{equation*}
\noindent Consequently,
 \begin{equation} \label{eq:bd}
     \|h_{T}\|_2 \leq \frac{ 2 \epsilon \sqrt{1+(\rho k-1)\mu}+2\mu\sqrt{\rho k}e }{(1+\mu+w \mu\sqrt{\rho k}-\mu \rho k(2 w \sqrt{\alpha}+1))}.
 \end{equation}
%by taking denominator positive, which is equivalent to 
The conditions in (\ref{boundK_new}) imply that $C_0$, $C_1$ and the denominator on the right hand side of (\ref{eq:bd})  are positive. 
\qed.
 \par An investigation into the choices of $T$, $\alpha$, $\rho$ and $w$ that result in smaller values for the RHS of (\ref{error bound T}) is presented in the following section.

%\section{Analysis of sufficient condition on sparsity} 
\begin{remark}
It is clear that the bound on sparsity $k$ in (\ref{boundK_new})  becomes less restrictive for small values of $\rho$ and $\mu$. For $\alpha=0$, the bound reduces to $\Bigg ( \frac{1}{2\sqrt{\rho}}\bigg(w+\sqrt{w^2+4(1+\frac{1}{\mu}}\bigg)\Bigg)^2$,  which is an increasing function of $w$. Similarly  for $\alpha \neq 0$, it is easy to verify that the $k$-bound is a decreasing function of $w$ and the largest $k$ bound is obtained at $\alpha=0$ and $w=1$.
\end{remark}
\begin{remark}
In the subsequent part, we refer to the bounds in (\ref{FriedBound}), (\ref{IETBound}),(\ref{ric_roc_Bound}) and (\ref{block_Bound}) as `global' bounds while the one in (\ref{error bound T}) as a `local'  bound. 
\end{remark}
\section{Analysis and comparison of `local` and `global' error bounds} \label{sec:local_bd}

In this section, we analyze the behaviour of the bound provided in (\ref{error bound T}) in terms of the associated parameters $\rho,\alpha$ and $w$. 

\par
\begin{figure}
    \centering
    \subfloat{\includegraphics[width=7cm,height=6cm]{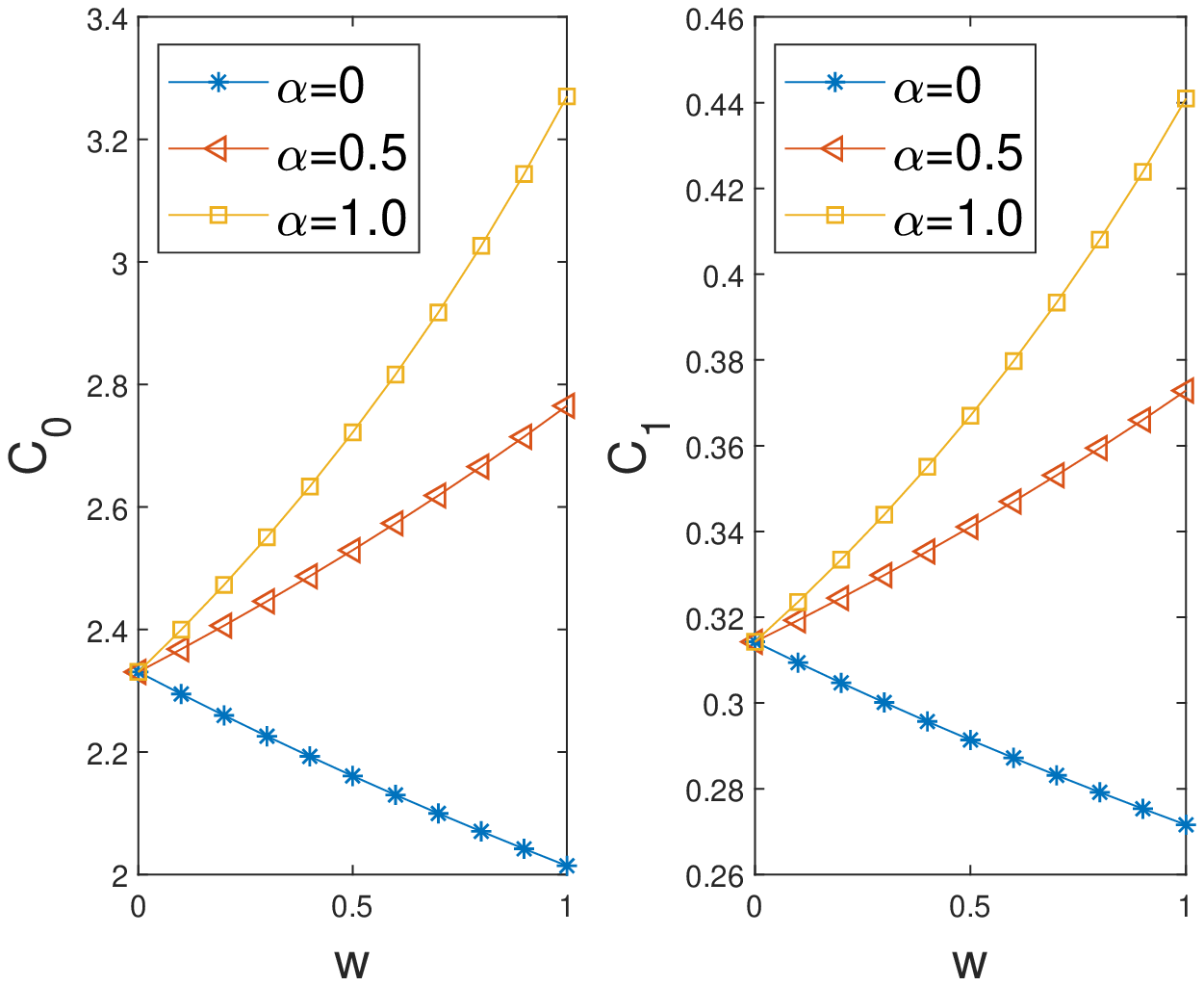} }
    \qquad
    \subfloat{\includegraphics[width=7cm,height=6cm]{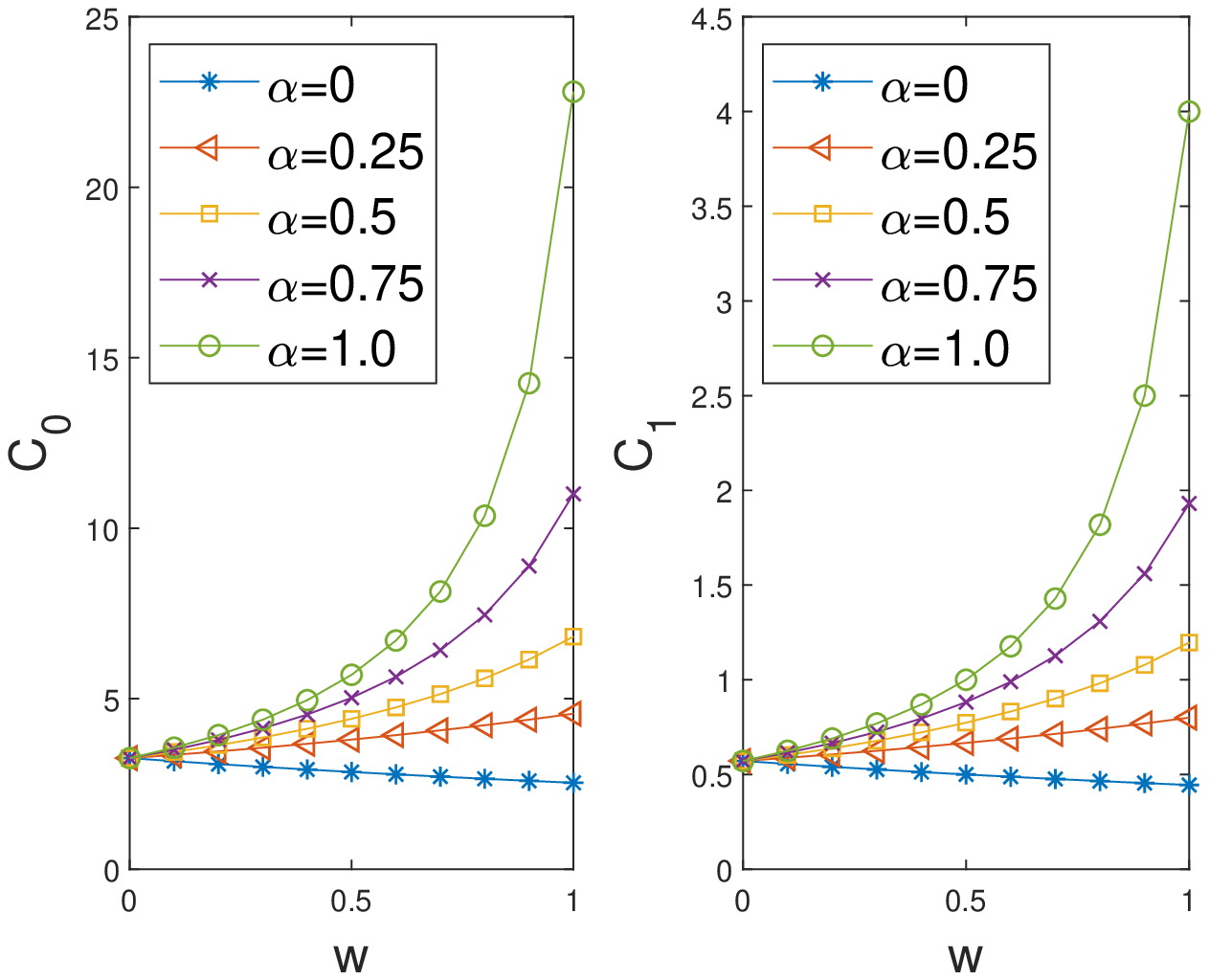} }
    \caption{The values of $C_0$ and $C_1$ for different $w$. The plots on the left and right panels are respectively for $\rho$=0.5 and $\rho$=1. For $\alpha=0$, the coefficients decrease with $w$ and for any $\alpha\neq 0$ increase with $w$. The smallest values of $C_0$ and $C_1$ are obtained around $\alpha=0,w=1$. Note that coefficients are smaller for smaller $\rho$.}
    \label{c0c1}
\end{figure}

 It may be noted that $\rho$ determines the relative size of $T$ with respect to the size of support of best $k$-term approximation of $x$. From (\ref{coeff}), it can be seen that the  coefficients $C_0$ and $C_1$ decrease with decrease in $\rho$. Again, when $\alpha=0$, the $w$-term in the denominator of coefficients is $w\mu\sqrt{\rho k}$, which being positive increases with $w$, making the coefficients decrease with $w$. When $\alpha\neq 0$, however, the $w$-term in the denominator of the coefficients is $w\mu\sqrt{\rho k}(1-\sqrt{\alpha\rho k})$, which is negative as $\alpha\rho k\geq 1$, since  $T$ is a nonempty subset of $[n]$. As a result, it  decreases with increase in $w$, which makes the coefficients increase with $w$. The stated behaviour of coefficients can be seen in Fig \ref{c0c1}. In generating the plots in this figure, as an example, we have taken $k=4$, considering the coherence of the underlying matrix as $\mu=0.1$,  like in \cite{liu2014compressed}. Note that $k=4$ and $\rho=0.5$ lead to three possibilities for the values of $\alpha$, viz, $0,0.5$ and $1$. Similarly, when $\rho=1$, $\alpha$ takes $0,0.25,0.5,0.75$ and $1$ as possible values. 
\par Since $C_0\epsilon$ term is controlled by $\epsilon$, for the reconstruction error to be small in (\ref{error bound T}), the multiplier $e$ of $C_1$ should be small, where  $e=w\|x-x_k\|_1 +(1-w)\|x_{T^c \cap T_0^c}\|_1+\|x_{T^c \cap T_0}\|_1$. The first two terms are small in the error expression as $x_k$ is the best-$k$-term approximation of $x$ and $T_0$ is the support of $x_k$ which can be further reduced by choosing optimal $w$. In order to make $e$ small, we need $\|x_{T^c \cap T_0}\|_1$ to be small. This is possible if $x_T$ contains the largest components of $x_k$ or the cardinality of $T^c \cap T_0$ is as small as possible. The latter case can happen when $\alpha\rho$ is close to 1. 
In application like in interior tomography, this condition translates to interior portion possessing dominating pixels. As the objective of the paper not is related to tomography, we do not go into the details of it any further. %In practical terms, if ROI contains dominating pixels of the image, then $e$ becomes small.

\par
\begin{figure}
    \centering
    \subfloat{\includegraphics[width=7cm,height=6cm]{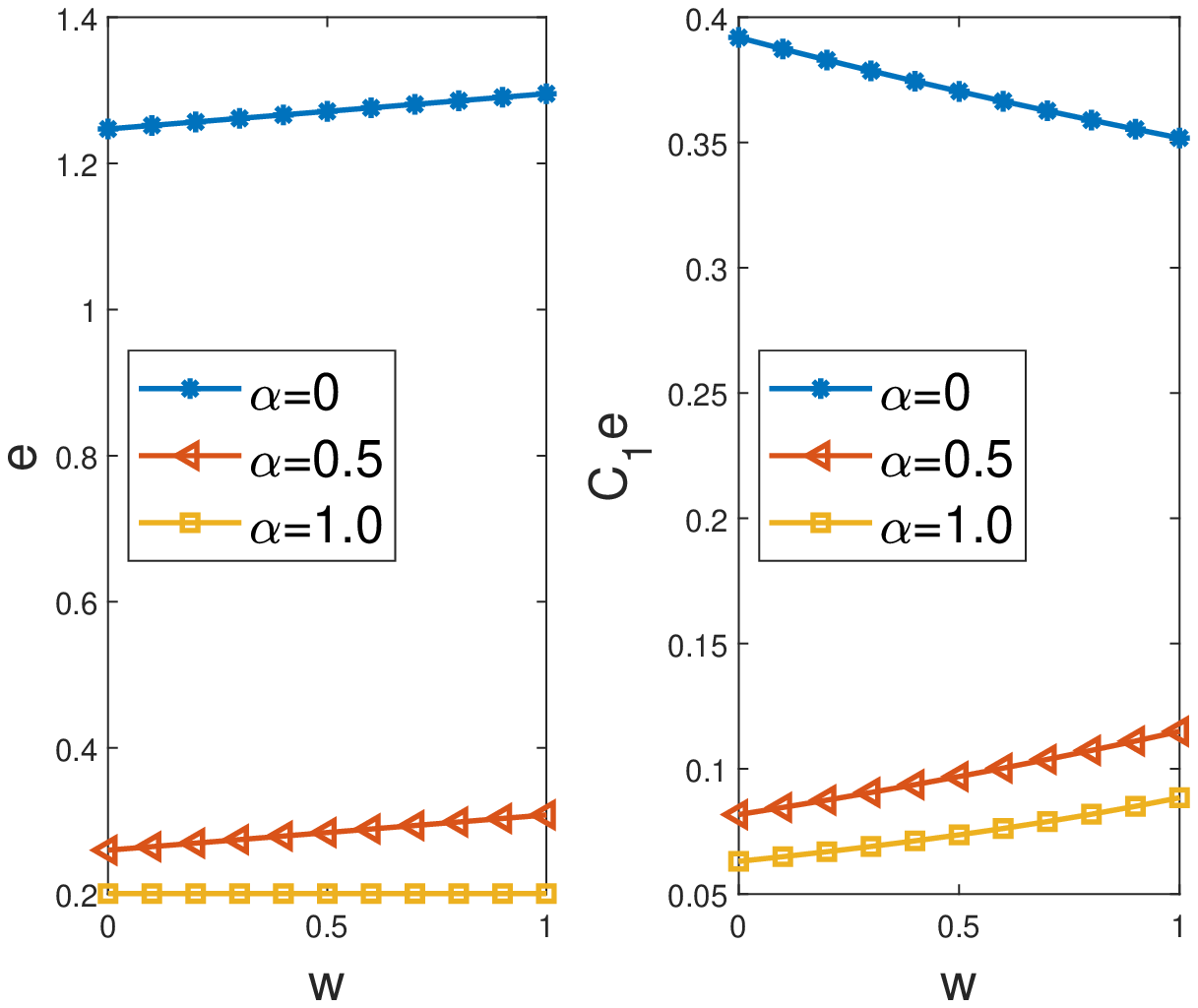}}
    \qquad
    \subfloat{\includegraphics[width=7cm,height=6cm]{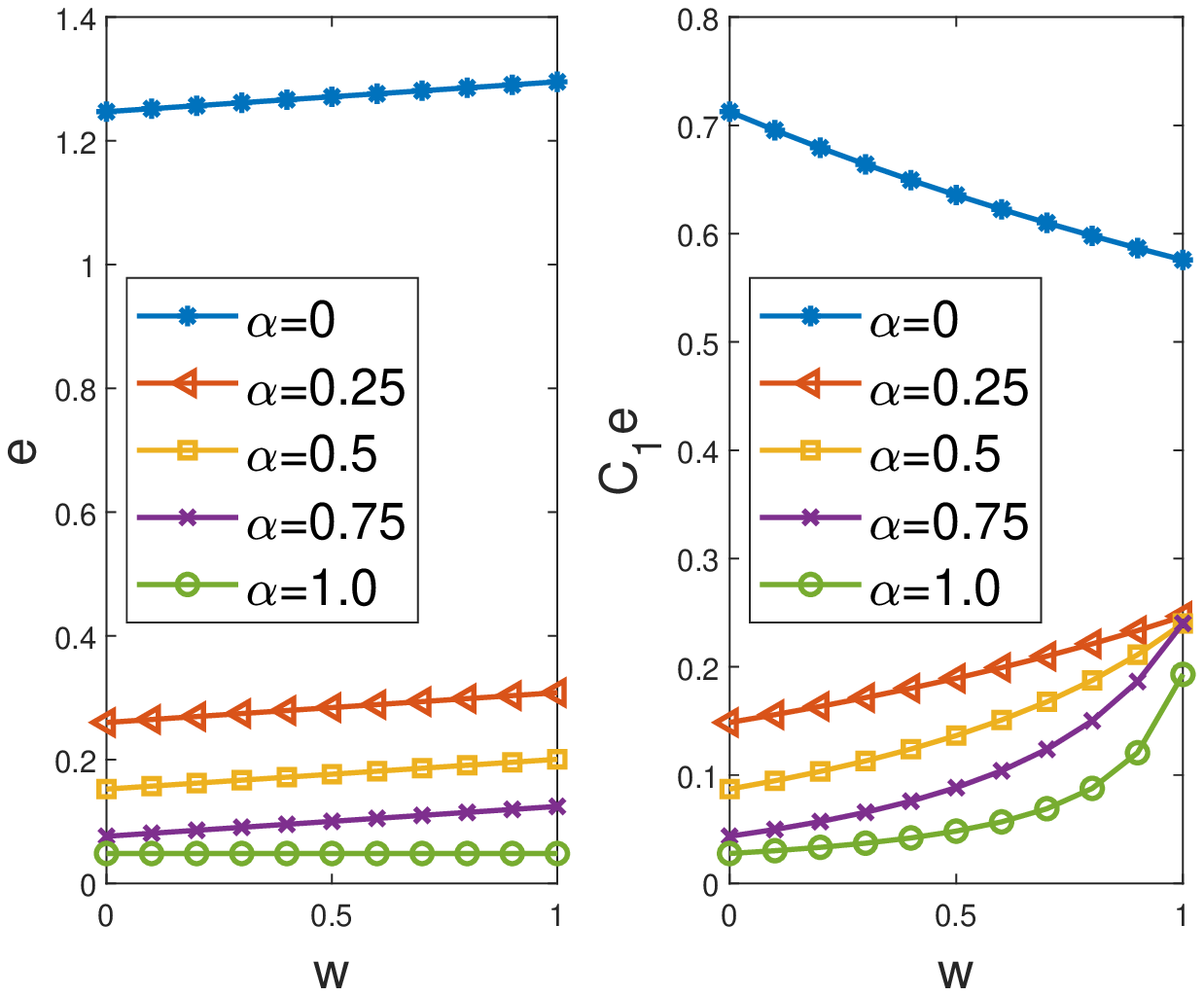}}
    \caption{The values of $e$ and $C_1e$  for different $w$ and $\alpha$ values. These plots indicate that the dominant error term in (\ref{error bound T}) takes smaller values when $\alpha=1$ and $w =0$.}
    \label{e_plot}
\end{figure}

 We have computed the values of $e$ and $C_1 e$ against different $\alpha$ and $w$, which are depicted in Fig. \ref{e_plot}. Here $x$ has been taken to be a normalized vector from Gaussian distribution with previously stated values for $k$ and $\mu$. 
%which is 5-sparse, $k=4$ and $\mu=0.1$({\color{red}{doubt}}). 
The plots in this figure  indicate that $e$ is large for $\alpha=0$, and for a given $\alpha$, it increases with $w$. The smallest value for $e$ is obtained at $\alpha=1$. Again $C_1 e$ is large for $\alpha=0$. For smaller values of $\alpha$, $C_1 e$ decreases with $w$ and for $\alpha=0$, $C_1e$ decreases with $w$. This is due to the effect of $C_1$. Overall, the error-bound takes its minimum value around $\alpha=1$ and $w=0$. 

\subsection{Comparison of bounds on $k$}
\begin{figure}
    \centering
    \subfloat{\includegraphics[width=7cm,height=7cm]{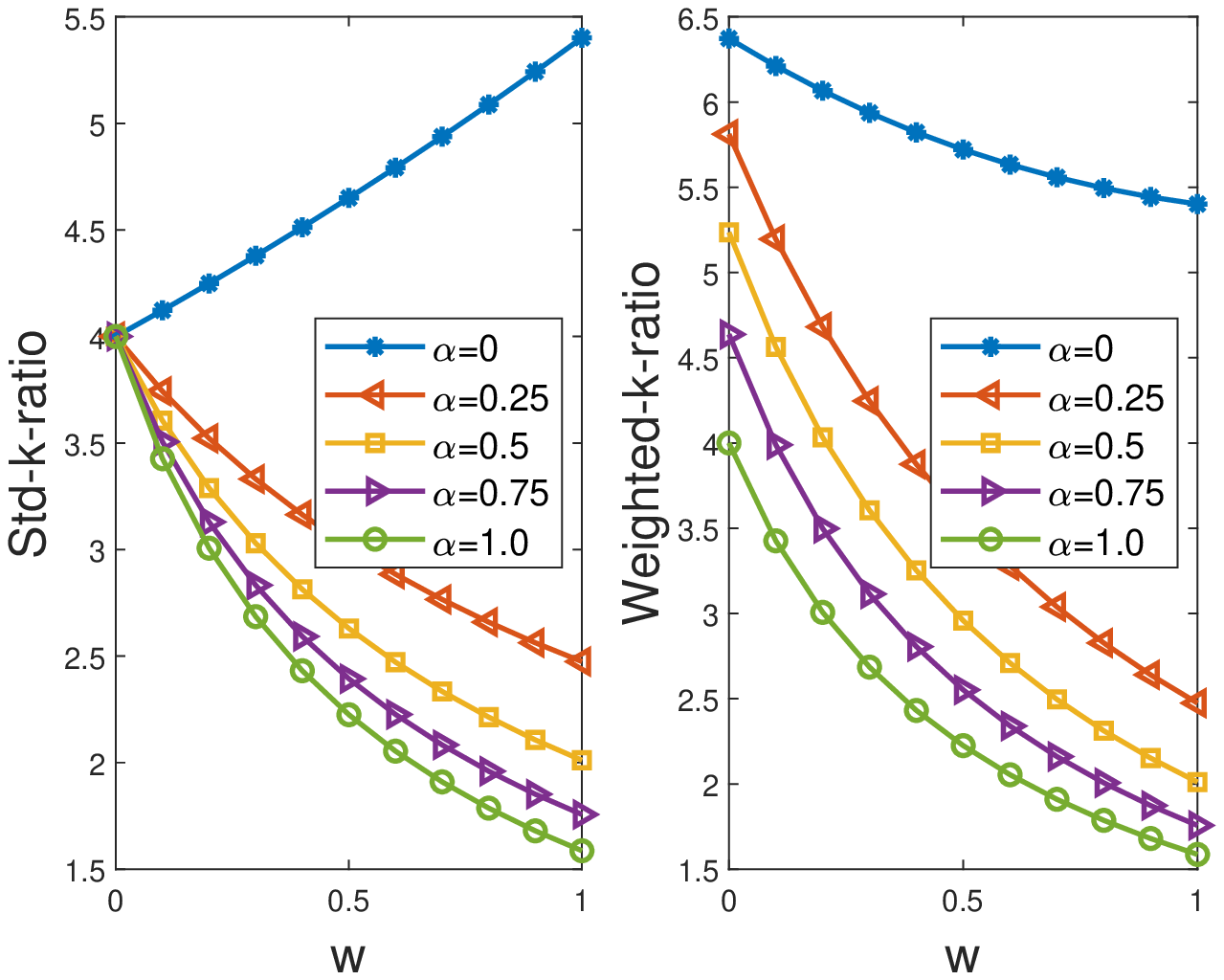}}
    \qquad
    \subfloat{\includegraphics[width=7cm,height=7cm]{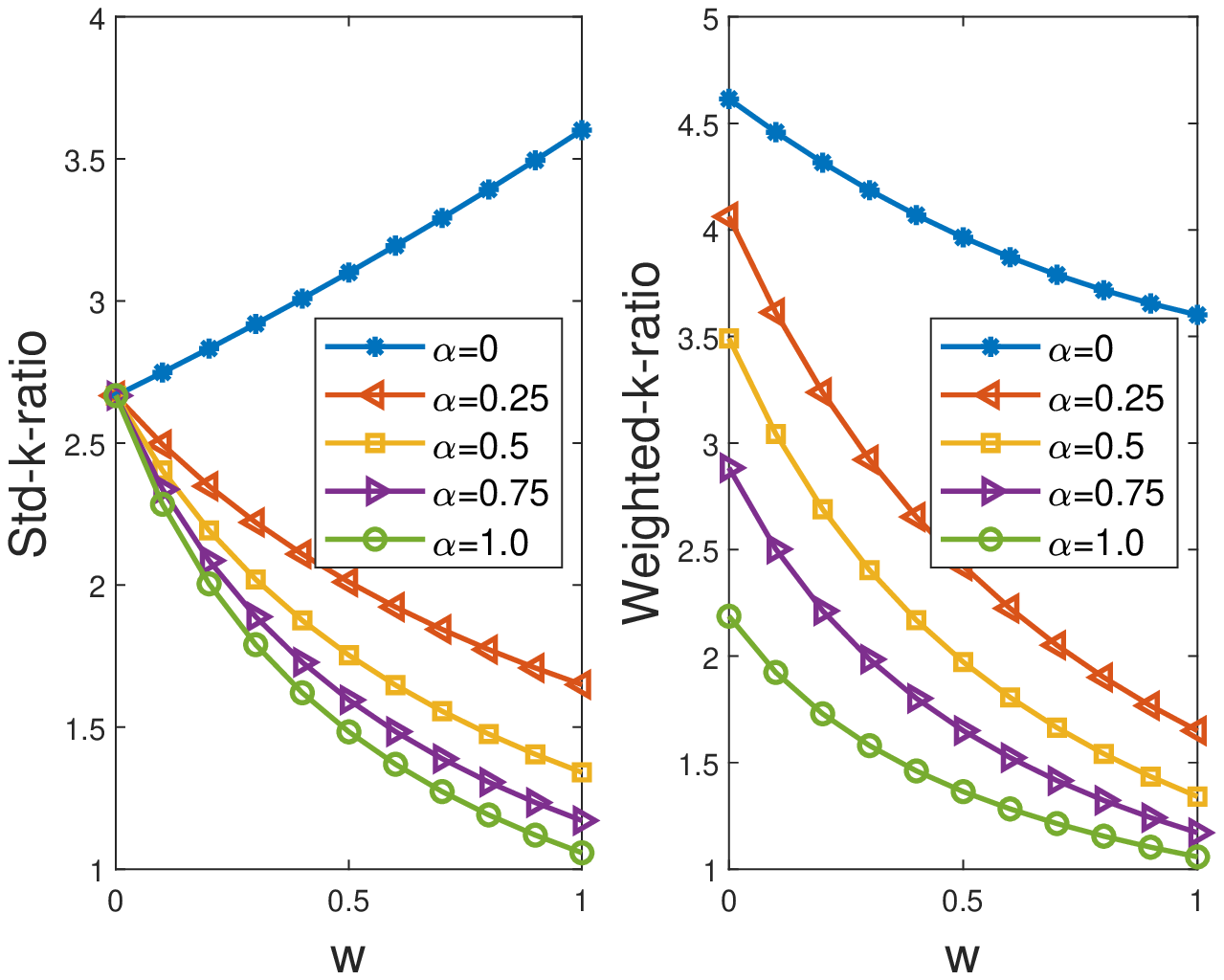}}
    \caption{Comparison of $k$-ratios for $\rho=0.5$ and $\rho=0.75$ respectively. Note that in all cases $k$-ratio is strictly greater than 1, which indicates that our sufficient condition on sparsity $(k)$ is much weaker than that of the standard and weighted one norm cases.}
    \label{k_plot}
\end{figure}

 Though the objective of present work is to propose an error bound restricted to a nonempty subset $T\subset[n]$, a natural question arises whether our $k$ bound in (\ref{boundK_new}) is better behaved than the ones associated with the global error bounds in the standard case (\ref{boundK_classical}) and in the weighted case (\ref{boundK_iet}). We compare the $k$-bounds  numerically in terms of a $k$-ratio. By a $k$-ratio in the standard case we  mean the right hand side of our $k$-bound in (\ref{boundK_new}) 
 divided by its counterpart in (\ref{boundK_classical}).  Similarly we  consider the  k-ratio in the weighted one norm case with respect to (\ref{boundK_iet}).
  Fig. \ref{k_plot} provides a comparison of $k$-ratios, which has been generated by considering, as examples,  $\mu=0.1$, $\rho=0.5$ and $\rho=0.75$.
 It can be seen in both the cases that the $k$-ratios are strictly greater than 1 for all $w\in[0,1]$ and for all $\alpha\in[0,1]$. It is clear from the graphs that,  even for larger values of $\rho$, our $k$ bound becomes less pessimistic  than that of  (\ref{boundK_classical})  and (\ref{boundK_iet}) for smaller values of $w$ for all $\alpha$. But since we are interested in finding error bounds on a subset $T$, which has a smaller size than the support of the best $k$-term approximation $T_0$, we do not consider the cases for larger values for $\rho$. 

\subsection{Comparison of bounds on error}
The local error bound in Theorem \ref{theorem_hT}  (that is, (\ref{error bound T})) becomes relevant if its right hand side is less than that of the global bounds in  Theorems \ref{Fried}, \ref{IETbound}, \ref{ric_roc} and \ref{block_thm}.
Since the stated right hand sides are functions of the solution vector $x$ with various associated parameters and different underlying conditions, comparing them for a general solution vector $x\in \mathbf{R}^n$  does not look practical. In view of this, we try to compare the coefficients (i.e, $C_0$ and $C_1$ with their respective counterparts in Theorems \ref{Fried}, \ref{IETbound}, \ref{ric_roc} and \ref{block_thm}) when the associated error expressions coincide. It is clear that the associated error parts coincide if $\|x_{T^c \cap T_0}\|_1$ is small. That is, $T$ contains largest magnitude entries from $T_0$. Consider a particular case for this in which $T=T_0$. In such a case, both the error terms associated with  coefficients in global as well as local error bounds coincide. Hence in this case it is enough to compare the corresponding coefficients $C_0$ and $C_1$ (that is, $C_0$ and $C_1$ are respectively compared to $C_0^{(i)}$ and $C_0^{(i)}$ for $i=1,2,3,4$).
\par
\begin{figure}
    \centering
    \subfloat{\includegraphics[width=7.5cm,height=6cm]{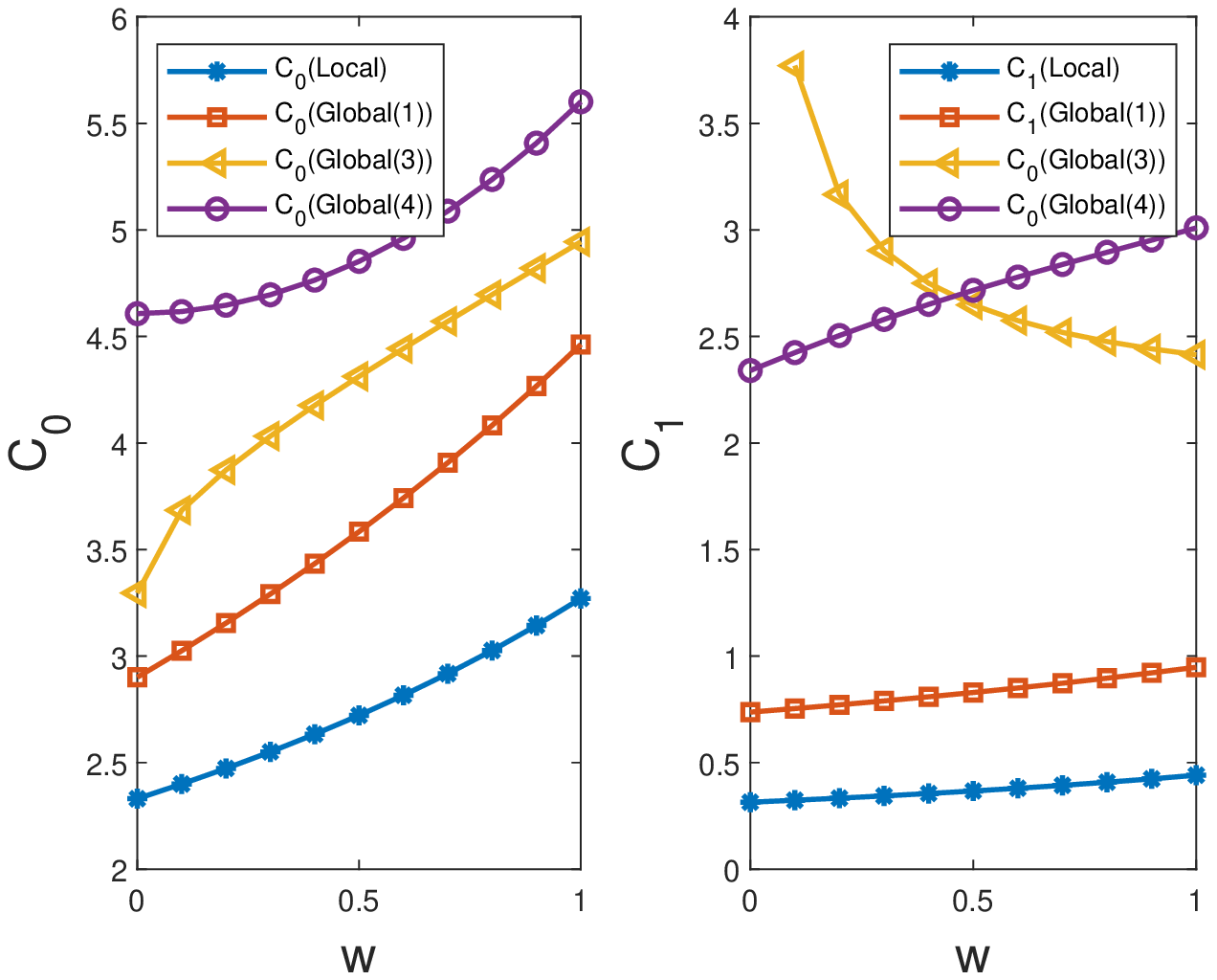}}
    \qquad
    \subfloat{\includegraphics[width=7.5cm,height=6cm]{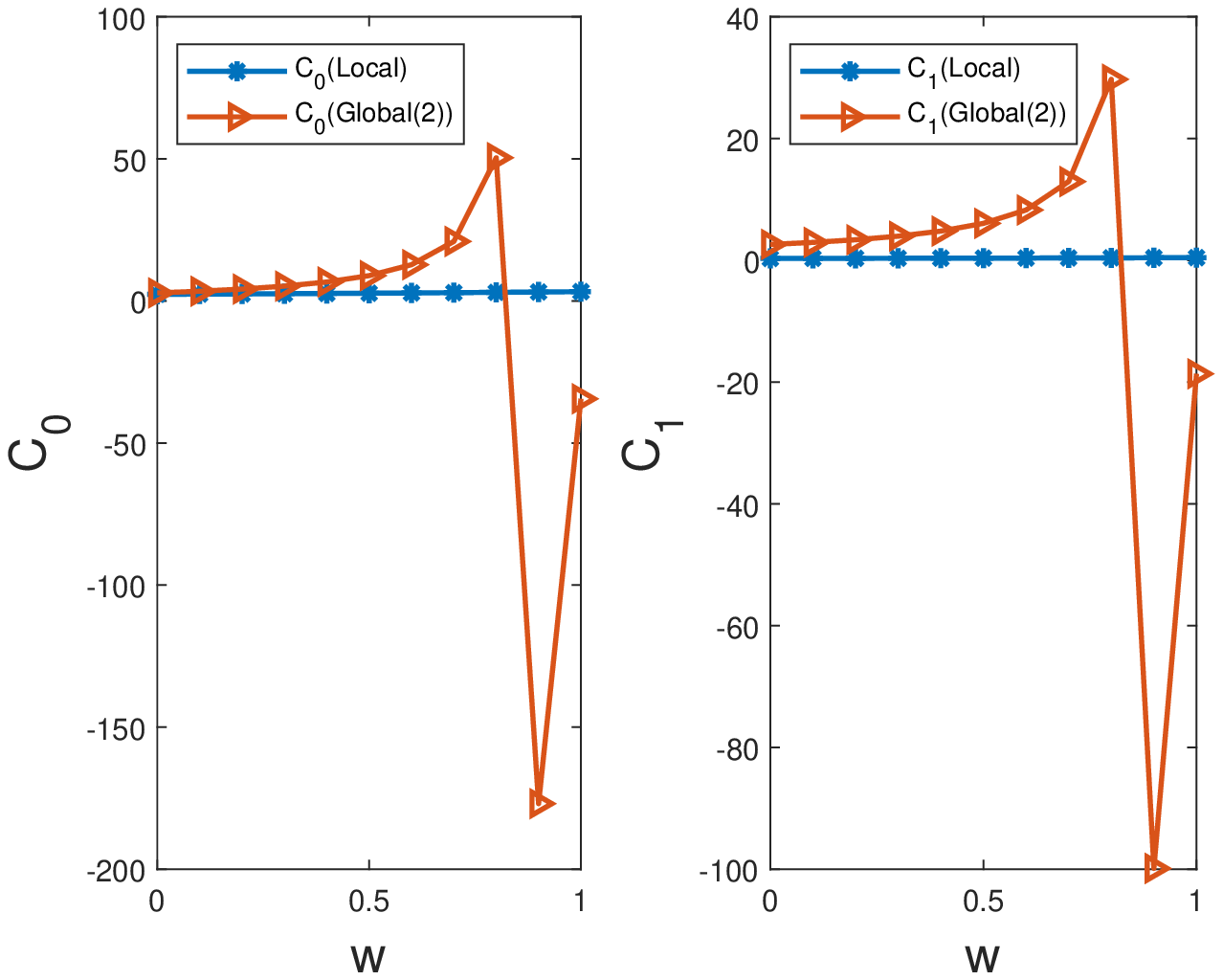}}
   \caption{Comparison of the coefficients when the terms associated with them coincide. Note that both coefficients corresponding to our local bound are less than the corresponding ones in other bounds, when the partial support $T$ contains the indices corresponding to largest magnitude entries of the solution $x$. }
    \label{comp_plot}
\end{figure}

 The coefficients in (\ref{coeff_iet}) and (\ref{coeff}) can be compared directly as both of them are in terms of mutual coherence parameter $\mu$. The coefficients in (\ref{coeff_frd}), (\ref{coeff_ric_roc}) and (\ref{coeff_block}) are in terms of the RIC and ROC. In order to compare them with the corresponding ones in (\ref{coeff}) which are in terms of $\mu$, we use the upper bounds \cite{candes2005decoding}\cite{elad2010sparse}: $\theta_{k,k}\leq \delta_{2k}$  and $\delta_k\leq (k-1)\mu$. For comparing with (\ref{coeff_frd}) we need a constant $a$ such that $a>1$ and $a k\in\mathbb{Z}$. We take a simple choice $a=2$. Again, for comparison with (\ref{coeff_ric_roc}), we need two constants $a$ and $b$. We set these to  $k$, which is permissible. Finally, for comparing with (\ref{coeff_block}), we need a constant $t>d$. In our case, as $d=1$, we take a simple choice $t=2$. The plots in Fig.  \ref{comp_plot}, comparing the coefficients in the stated setting, are denoted through the legends `Local', `Global(1), Global(2), Global(3) and Global(4), which stand respectively for the coefficients in  (\ref{coeff}), (\ref{coeff_iet}),(\ref{coeff_frd}),(\ref{coeff_ric_roc}) and (\ref{coeff_block}). 
From  Fig \ref{comp_plot},  it is clear that our coefficients being much small imply that the right hand side of our bound is small compared to their global bounds at least in the case where $T$ contains the indices corresponding to largest magnitude entries in $x$ for all $w\in[0,1]$. As highlighted already, comparison of bounds in other cases does not look possible. 

In generating Fig \ref{comp_plot}, we have taken $\mu=0.1$ and $k=2$. 
It may be noted that the coefficients in (\ref{coeff_frd}) are negative for $w>0.8$ and in (\ref{coeff_ric_roc}), $w=0$ gives $s=0$ which makes $C_0^{(3)}$ in (\ref{coeff_ric_roc}) not defined at $w=0$.

\section{Conclusion}
The present work has proposed a local recovery bound for prior support constrained compressed sensing, while the existing bounds are global in nature. In particular, an error estimate restricted to prior support, providing recovery guarantee, has been provided along with a bound on the sparsity of the solution to be recovered.  \vspace{3mm}

{\bf Acknowledgments}:\\
The first author is thankful to the UGC, Government of India, (JRF/2016/409284) for its financial support. The second author gratefully acknowledges the support received from the MHRD, Government of India.
\normalem
\bibliographystyle{plain}
\bibliography{bibliography.bib}

\end{document}